\documentclass[11pt]{article}
\usepackage{graphicx}
\usepackage{subfigure}
\usepackage{tikz}
\usepackage{tikz-network}
\usepackage{float}
\usetikzlibrary{positioning,petri}
\usepackage{amsfonts}
\usepackage{amsbsy}
\usepackage{amssymb}
\usepackage{amsmath,amsthm}
\usepackage{verbatim}
\usepackage{amsfonts}
\usepackage{parskip}
\usepackage[T1]{fontenc}
\usepackage{multicol}
\usepackage[onehalfspacing]{setspace}
\usepackage{color}
\usepackage[autostyle]{csquotes}
\usepackage{pgf}
\usepackage{mathtools}
\usepackage{hyperref}
\usepackage{thmtools}
\usepackage{caption}
 
\declaretheoremstyle[%
  spaceabove=-6pt,%
  spacebelow=6pt,%
  headfont=\normalfont\itshape,%
  postheadspace=1em,%
  qed=\qedsymbol%
]{mystyle} 
\declaretheorem[name={Proof},style=mystyle,unnumbered,
]{prf}
\hypersetup{
    colorlinks=true,
    linkcolor=blue,
    filecolor=blue,      
    urlcolor=blue,
    citecolor=blue,
    pdftitle={Overleaf Example},
    pdfpagemode=FullScreen,
    }

\newtheorem{thm}{Theorem}[section]
\newtheorem{prop}{Proposition}[section]

\newtheorem{lem}[thm]{Lemma}

\newtheorem{rem}[prop]{Remark}
\newtheorem{exm}[thm]{Example}

\RequirePackage{pgfkeys}
\pgfkeys{/ytableau/options/.is family}

\pgfkeys{/ytableau/options, boxsize/.value required,  boxsize/.code = {%
\pgfkeysalso{nosmalltableaux}%
 \compare@yT{#1}{normal}%
 \ifeq@yT
 \xdef\macro@boxdim@yT{\expandonce@yT\boxdim@normal@yT}%
 \else
 \xdef\macro@boxdim@yT{#1}%
 \fi
 } }
 \pgfkeys{/ytableau/options, textmode/.value forbidden,
 textmode/.code = \set@textmode@yT,
mathmode/.value forbidden,
mathmode/.code = \set@mathmode@yT,
 }
 
\def\set@mathmode@yT{
 \gdef\skipin@yT{$}
 \gdef\skipout@yT{$}
 \def\smallfont@yT{\scriptstyle} } 
 \set@mathmode@yT

\makeatletter
\let\@fnsymbol\@arabic
\makeatother

\title{Some invariants related to threshold and chain graphs}
\author{Rameez Raja$^{1}$ and Samir Ahmad Wagay$^{2}$\footnote{$^{, 2}$Department of Mathematics, National Institute of Technology Srinagar, Jammu and Kashmir, India. Email: rameeznaqash@nitsri.ac.in, samir\_03phd20@nitsri.net,
Corresponding author:$^1$}}

\begin{document}
\maketitle
\begin{abstract} 
 Let $G = (V, E)$ be a finite simple connected graph. We say a graph $G$ \textit{realizes} a code of the type $0^{s_1}1^{t_1}0^{s_2}1^{t_2}\cdots0^{s_k}1^{t_k}$ if and only if $G$ can obtained from the code by some rule. Some classes of graphs such as \textit{threshold} and \textit{chain} graphs realizes a code of the above mentioned type. In this paper, we develop some computationally feasible methods to determine some interesting graph theoretical invariants. We present an efficient algorithm to determine the \textit{metric dimension} of threshold and chain graphs. We compute \textit{threshold dimension} and \textit{restricted threshold dimension} of threshold graphs. We discuss \textit{$L(2, 1)$-coloring} of threshold and chain graphs. In fact, for every threshold graph $G$, we establish a formula by which we can obtain the \textit{$\lambda$-chromatic number} of $G$. Finally, we  provide an algorithm to compute the \textit{$\lambda$-chromatic number} of chain graphs.  
  \end{abstract}

2020 AMS Classification Code: 05C12, 05C15, 05C85. 
 
\section{Introduction}
This paper is concerned with \textit{threshold} and \textit{chain} graphs. Recall that a graph $G$ is called a \textit{chain} graph if it is bipartite and the neighbourhoods of vertices in each partite set form a chain with respect to set inclusion. The color classes $U$, $V$ of a chain graph $G$ can be partitioned into $h$ non-empty cells $U_1, \cdots , U_h$ and $V_1, \cdots , V_h$ such that $ N(u) = V_1 \cup \cdots \cup V_{h - i+1}$ for any vertex $u\in U_i$ with $1\leq i\leq h$, $N(u)$ denotes the neighbourhood of a vertex $u$ of $G$. Equivalently, chain graphs are $\{2K_2, C_3, C_5\}$-free graphs. Chain graphs present a broad class of graphical models for description of conditional independence structures, including computer networks, social networks, communication networks, information networks, software designs and Bayesian networks.

On the other hand, \textit{threshold} graphs are $\{P_4, C_4, 2K_2\}$-free graphs. The concept of threshold graphs was first established by the authors Chvátal and Hammer \cite{CH} and Henderson and Zalcstein \cite{HZ}. Threshold graphs play an important role in graph theory as well as in several applied areas such as psychology, computer science, scheduling theory, etc. These graphs can be used to control the flow of information between processors, much like the traffic lights used in controlling the flow of the traffic. These graphs can be viewed as special cases of some wider classes of graphs like cographs, split graphs and interval graphs.

It is important to mention here that  threshold and chain graphs can be generated by some binary code. In fact, given a binary code of the type $0^{s_1}1^{t_1} \cdots0^{s_k}1^{t_k}$, where $s_i$ and $t_i$ are positive integers with $1\leq i \leq k$, we can construct these graphs. A formal definition is presented in section 2. In this paper, we are mainly interested to determine some graph theoretical invariants such as \textit{metric dimension}, \textit{threshold dimension}, \textit{restricted threshold dimension} and \textit{$\lambda$-chromatic number} from a binary generating code  $0^{s_1}1^{t_1} \cdots0^{s_k}1^{t_k}$. We focus on developing some computationally feasible methods to determine the invariants discussed in the succeeding sections. 

We have organized the theory and the results involved in this investigation as follows. In section $2$, we present some preliminary information related to threshold and chain graphs. We discuss some results regarding \textit{metric dimension} of threshold and chain graphs. We present an algorithm for computing the \textit{metric dimension} of a threshold graph and a chain graph from their respective binary generating code. In section $3$, we provide bounds for the \textit{threshold dimension} and describe the method (Algorithm-(a)) for computing the \textit{restricted threshold dimension} of a threshold graph. Finally, we discuss the \textit{$\lambda$-chromatic number} of a threshold graph and introduce an algorithm (Algorithm-(b)) to determine the \textit{$\lambda$-chromatic number} of chain graphs.

\section{Metric dimension}
A threshold graph can be defined in several different ways as can be seen in \cite{MP}. We follow the definition through binary generating codes, that is, a threshold graph $G(b)$ is obtained from a binary code of the form $b = (b_1b_2 \cdots b_n)$ in the following way,
\begin{itemize} 
 \item[(i)] Assume $G_1 = G(b_1) = K_1$,  a single isolated vertex.
 \item[(ii)] For $i\geq2$, let $G_{i-1}= G(b_1b_2...b_{i-1})$. Then the graph $G_i = G(b_1b_2...b_{i-1}b_i)$ is obtained by adding $b_i$ vertices to $G_{i-1}$ (already constructed). Note that $b_i = 0$ if the vertex $v_i$ is added as an isolated vertex to $G_{i-1}$, and $b_i = 1$ if $v_i$ is added as a dominating vertex to $G_{i-1}$. This representation is called a binary generating code for the threshold graph and is written as,
 $$b ~=~ (0^{s_1}1^{t_1})(0^{s_2}1^{t_2}) \cdots (0^{s_k}1^{t_k}), \mbox{where} ~s_i, t_i > 0 ~ \mbox{for} ~ 1\leq i\leq k.$$
 \end{itemize}
We restrict ourselves with connected graphs only, so $b_k= 1$. A simplified illustration of threshold graph is shown in Figure 1-(a).\\
Chain graphs are analogues to threshold graphs in the class of all connected bipartite graphs of specified order and size. By deleting all the edges that belong to the complete part of a threshold graph, we obtain the chain graph. Accordingly, a chain graph can also be generated by a binary code in  the following way,
\begin{itemize}
\item[(i)] Assume $G_1 = G(b_1) = K_1$,  a single isolated vertex.
 \item[(ii)] For $i\geq2$, let $G_{i-1}= G(b_1b_2...b_{i-1})$. Then the graph $G_i = G(b_1b_2...b_{i-1}b_i)$ is obtained by adding $b_i$ vertices to $G_{i-1}$ (already constructed). Note that $b_i = 0$ if the vertex $v_i$ is added as an isolated vertex to $G_{i-1}$ and $b_i = 1$ if $v_i$ is added as a dominating vertex to the only previously added isolated vertices of $G_{i-1}$. A representation of the chain graph is shown in Figure 1-(b).
 \end{itemize}

 \begin{figure}[H]
   \begin{tikzpicture}
  \Vertex[y=0, label=$U_5$, color=white, position=left]{A}\Vertex[y=1,label=$U_4$, color=white, position=left]{B}\Vertex[y=2,label=$U_3$, color=white, position=left]{C}\Vertex[y=3,label=$U_2$, color=white, position=left]{D}\Vertex[y=4,label=$U_1$, color=white, position=left]{E}\Vertex[x=3.5,y=0, label=$V_5$, color=black, position=right]{F}\Vertex[x=3.5,y=1,label=$V_4$, color=black, position=right]{G}\Vertex[x=3.5,y=2,label=$V_3$, color=black, position=right]{H}\Vertex[x=3.5,y=3,label=$V_2$, color=black, position=right]{I}\Vertex[x=3.5,y=4,label=$V_1$, color=black, position=right]{J}
  \Edge(F)(J)\Edge(E)(J)\Edge(E)(I)\Edge(E)(G)\Edge(E)(F)\Edge(E)(H)\Edge(D)(I)\Edge(D)(H)\Edge(D)(F)\Edge(C)(H)\Edge(G)(E)\Edge(G)(D)\Edge(G)(C)\Edge(G)(B)\Edge(F)(A)\Edge(F)(B)\Edge(F)(C)\Edge(F)(D)\Edge(F)(E)\Edge[bend=35](J)(H)\Edge[bend=55](J)(G)\Edge[bend=55](J)(F)\Edge[bend=65](I)(F)\Edge[bend=70](I)(G)\Edge[bend=70](H)(F)
\end{tikzpicture}
\centering{(a)}
\hfill
\begin{tikzpicture}
  \Vertex[y=0, label=$U_5$, color=white, position=left]{A}\Vertex[y=1,label=$U_4$, color=white, position=left]{B}\Vertex[y=2,label=$U_3$, color=white, position=left]{c}\Vertex[y=3,label=$U_2$, color=white, position=left]{D}\Vertex[y=4,label=$U_1$, color=white, position=left]{E}\Vertex[x=3.5,y=0, label=$V_5$, color=black, position=right]{F}\Vertex[x=3.5,y=1,label=$V_4$, color=black, position=right]{G}\Vertex[x=3.5,y=2,label=$V_3$, color=black, position=right]{H}\Vertex[x=3.5,y=3,label=$V_2$, color=black, position=right]{I}\Vertex[x=3.5,y=4,label=$V_1$, color=black, position=right]{J}
  \Edge(E)(J)\Edge(E)(I)\Edge(E)(G)\Edge(E)(F)\Edge(E)(H)\Edge(D)(I)\Edge(D)(H)\Edge(D)(F)\Edge(C)(H)\Edge(G)(E)\Edge(G)(D)\Edge(G)(C)\Edge(G)(B)\Edge(F)(A)\Edge(F)(B)\Edge(F)(C)\Edge(F)(D)\Edge(F)(E)
\end{tikzpicture}
{(b)}
\caption{Representation of threshold and chain graphs generated by a  code  $0101010101$}
\end{figure}
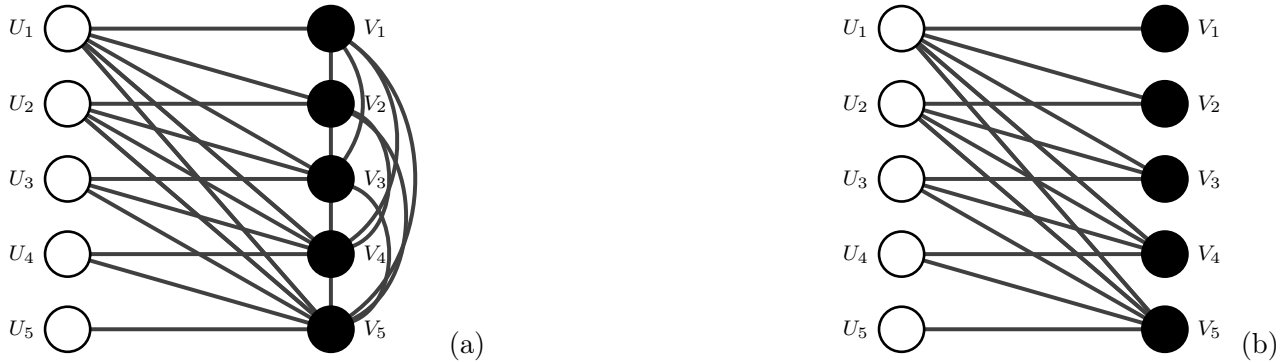
Note that, we call a part $0^{s}1^{t}$ of a binary generating code as a \textit{string}. The power bit $0^s$ represents the independent part (white class) and $1^t$ represents the complete part (black class) of vertices of the graph $G$ realized by $0^{s_1}1^{t_1}0^{s_2}1^{t_2}\cdots0^{s_k}1^{t_k}$.

Let $\mathcal{W}  = \{v_1, v_2, \cdots, v_k\}$ be an ordered set of vertices of $G$ and let $w\in V(G)$. The representation of $w$ with respect to $\mathcal{W}$ is the ordered $k$-tuple $d(w|\mathcal{W})=\big(d(w,v_{1}),d(w,v_{2}), \cdots, d(w,v_{k})\big)$, where $d(w,v_i)$ represents distance between the vertex $w$ and vertices $v_i$, where $1\leq i \leq k$. If every two vertices of $G$ have distinct representations with respect to $\mathcal{W}$, then the set $\mathcal{W}$ is called a \textit{resolving set} (locating set) for $G$. The \textit{metric dimension} of $G$ is the minimum cardinality of a resolving set of $G$ denoted by $\beta(G)$. A resolving set containing a minimum number of vertices is called a \textit{basis} (reference set) for $G$. 

 Resolving sets for graphs were introduced independently by Slater \cite{P} and Harary and Melter\cite{FR}, while the concept of a resolving set and that of metric dimension were defined much earlier in the more general context of metric spaces by Blumenthal in his monograph theory and applications of distance geometry. Graphs are special examples of metric spaces with their intrinsic path metric. The metric dimension has appeared in various applications in graph theory as diverse as pharmaceutical chemistry, combinatorial optimization \cite{AE}, mastermind game tactics \cite{V}, network discovery and verification \cite{ZFT}, robot navigation \cite{SBA} and so on. According to the authors \cite{LP}, the problem of determining the metric dimension of an arbitrary graph is NP-complete. The following interesting results have appeared in \cite{GLMO}.

\begin{thm}
If $G$ is a connected graph of order $n\geq 2$ and diameter $d$,  then $f(n, d)\leq \beta(G)\leq n-d$, where $f(n, d)$ is the least positive integer $k$ for
which $k + d^k\geq n$.
\end{thm}
\begin{thm}
Let $G$ be a connected graph of order $n\geq 2$. Then the following hold,
\begin{itemize}
\item[(i)] $\beta(G)=1$ if and only if $G$ is a path and,
\item[(ii)] $\beta(G)=n-1$ if and only if $G$ is a complete graph.
\end{itemize}
\end{thm}
In the following result, we determine the metric dimension of a threshold graph realized by a string of type $0^{s}1^{t}$.
  
 \begin{lem} If $G$ is a graph of order $s+t$ realized by the string $0^{s}1^{t}$, then $\beta(G) =  s+t-2$, where $s >1$. 
 \begin{prf}
 We partition the vertex of $G$ as, $V(G) = U\cup V$, where $U~ = \lbrace u_1,u_2, \cdots,u_s\rbrace$ and $V = \lbrace v_1,v_2, \cdots, v_t\rbrace$  represents the independent and complete parts of $G$. Since $G$ is threshold, therefore $d(u_i, u_k) = 2$, for all $i \neq k$ and $d(u_i, v_j)= d(v_j, v_l)= 1$, for all  $j \neq l $. It follows that there are only two vertices in $G$ (one from each $U$ and $V$) whose distances are uniquely determined by the vertex set of $G$. Therefore the resolving set of $G$ contains exactly $s-1$ vertices of $U$ and $t-1$ vertices of $V$, that is, $\mathcal{W} = \lbrace u_1, u_2, \cdots, u_{s-1}, v_1, v_2, \cdots, v_{t-1}\rbrace$ contains all but two vertices of $G$ and the result follows.  
 \end{prf}
 \label{1}
  \end{lem}
\begin{prop}
Let $G_1$ be a graph realized by $0^{s_1}1^{t_1}$ with $\beta(G_1)= \alpha$ and let $G_2$ be a graph realized by $0^{s_1}1^{t_1}0^{s_2}1^{t_2}$. Then the metric dimension of $G_2$ is one of the following.
\begin{itemize}
\item[(i)] If  $s_2 > 1, t_2 > 1$ and $G_1$ is not a star graph, then $\beta(G_2) = \alpha + (s_2 + t_2 - 2)$, otherwise it is $\alpha + (s_2 + t_2 - 1)$.
\item[(ii)] If  $s_2 > 1, t_2 = 1$ and $G_1$ is not a star graph, then $\beta(G_2) = \alpha + (s_2 - 1)$, otherwise it is $\alpha ~+~s_2$.
\item[(iii)] If  $s_2 = 1, t_2 > 1$, then $\beta(G_2) = \alpha + t_2$.  
\item[(iv)] If $s_2 = 1, t_2 = 1$, then $\beta(G_2) = \alpha + 1$.  
\end{itemize}
 \begin{prf}
As before, we partition the vertex set of $G_2$ as, $V(G_2) = U\cup V$, where $U=\lbrace(x_1,x_2, \cdots, x_{s_1})\cup (z_1, z_2, \cdots, z_{s_2})\rbrace = \lbrace I_1\cup I_2\rbrace$ and $V = \lbrace(y_1, y_2, \cdots, y_{t_1})\cup(u_1, u_2, \cdots, u_{t_2})\rbrace = \lbrace J_1\cup J_2\rbrace$. The sets $I_1, I_2$ and $J_1, J_2$ respectively represents the independent and complete parts of $G_2$. For all $i, j, k$ and $l$, $d(x_i, y_j) = d(x_i, u_l)= d(z_k, u_l) = 1$ and $d(x_i, z_k) = d(y_j, z_k) =2$. 
 
 To prove (i), assume $G_1$ is not a star graph. Then the distance vector of all vertices of $I_1$ is of the form $(2, 2,\cdots,2, 1, 1,\cdots,1,~ 2, 2, \cdots,  2, 1, 1,\cdots,1)$, where as the distance vectors of all vertices of $J_1$, $I_2$ and $J_2$ are of the forms $(1, 1, \cdots, 1, 1, 1, \cdots, 1, 2, 2, \cdots, 2, 1,1, \cdots, 1)$, $(2, 2,  \cdots, 2, 2, 2, \cdots, 2\\
 , 2, 2, \cdots, 2, 1, 1, \cdots,1)$ and $(1, 1, \cdots, 1, 1, 1,\cdots,1, 1, 1,\cdots, 1, 1, 1, \cdots, 1)$ respectively. Therefore, $\mathcal{W}=\lbrace x_1, x_2, \cdots, x_{s_1-1}, y_1, y_2, \cdots, y_{t_1-1}, z_1,z_2,\cdots,z_{s_2-1}, u_1,u_2,\cdots,u_{t_2-1}\rbrace$ is the required resolving set and consequently, $\beta(G_2) = s_1-1+t_1-1+s_2-1+t_2-1 = \alpha + (s_2+t_2-2)$.\\
 $~~~~~~~$If $G_1$ is a star graph, then $V(G_2) =  \lbrace  x_1,x_2,\cdots,x_{s_1}\rbrace \cup \lbrace y_1 \rbrace \cup \lbrace z_1,z_2,\cdots,z_{s_2}\rbrace \cup \lbrace u_1,u_2,\cdots,u_{t_2}\rbrace$. By Lemma \eqref{1}, the resolving set of $G_1$ contain all but one vertices of $I_1$. Taking this into consideration, $d(z_k, x_{i}) = 2$ for all $i, k$ and $d(u_l, x_i ) = 1$ for all $i, l$. Thus, $\mathcal{W} = \lbrace x_1,x_2,\cdots,x_{s_1-1}, z_1,z_2
,\cdots,z_{s_2}, u_1\\
,u_2,\cdots,u_{t_2 - 1}\rbrace$ and therefore, $\beta(G_2) = s_1-1+s_2+t_2-1 = \alpha + (s_2+t_2-1)$.
 
 To prove (ii), we first assume that $G_1$ is not a star graph. Then by Lemma \eqref{1}, the resolving set of $G_1$ is $\lbrace x_1, x_2, \cdots, x_{s_1-1}, y_1, y_2, \cdots, y_{t_1-1}\rbrace$. By adding $s_2$ number of isolated vertices and a single dominating vertex $u_1$ to $G_1$, the new graph $G_2$ obtained from $G_1$ has a vertex set of the form  $\lbrace  x_1,x_2,\cdots,x_{s_{1}}\rbrace \cup \lbrace y_1, y_2, \cdots, y_{t_{1}} \rbrace \cup \lbrace z_1, z_2,\cdots, z_{s_{2}}\rbrace \cup \lbrace u_1 \rbrace$. However, $d(z_k,x_i) = d(z_k,y_j) =2 ~~\mbox{for all}~~i,j,k$ and the vertex $u_1$ is adjacent to every vertex of $G_2$. Therefore, the resolving set of $G_2$ contains the resolving set of $G_1$ and all but one vertices of $I_2$. Thus, $\mathcal{W} = \lbrace x_1,x_2,\cdots,x_{s_1-1}, y_1,y_2,\cdots,y_{t_1-1}\\
 , z_1,z_2,\cdots,z_{s_2-1}\rbrace$ and $\beta(G_2) = s_1-1+t_1-1+s_2-1 = \alpha + (s_2-1)$.\\  
 $~~~~~$ If $G_1$ is a star graph, then $V(G_2) =  \lbrace  x_1,x_2,\cdots,x_{s_1}\rbrace \cup \lbrace y_1 \rbrace \cup \lbrace z_1,z_2,\cdots,z_{s_2}\rbrace \cup \lbrace u_1 \rbrace$. For each $i$ and $k$, we have $d(z_k, x_i) =d(z_k, y_1)= 2$, and therefore, the resolving set contains all but one independent vertices of $G_2$. This implies, $\mathcal{W} = \lbrace x_1,x_2,\cdots,x_{s_1-1}, z_1,z_2,\cdots,z_{s_2}\rbrace$ and $\beta(G_2) = s_1-1+s_2 = \alpha + s_2$.
 
 For (iii), we consider the graph $G_2$ obtained from $G_1$ by adding a single isolated vertex $z_1$ and $t_2 $ number of dominating vertices. Then the vertex set of $G_2$ is of the form $\lbrace  x_1,x_2,\cdots,x_{s_1}\rbrace \cup \lbrace y_1,y_2,\cdots,y_{t_1}\rbrace \cup \lbrace z_1 \rbrace \cup \lbrace u_1,u_2,\cdots,u_{t_2}\rbrace$. For each $i, j$ and $l$, $d(z_1, x_i) =d(z_1, y_j) =2$ and $d(u_l, x_i) =d(u_l, y_j) =d(u_l,z_1) =1$. Therefore, the resolving set of $G_2$ contains the resolving set of $G_1$ and all the newly added dominating vertices of $G_2$. Thus, $\mathcal{W}=\lbrace x_1, \cdots, x_{s_1-1}, y_1, \cdots, y_{t_1-1}, u_1, u_2, \cdots, u_{t_2}\rbrace$  and $\beta(G_2) = s_1-1+t_1-1+t_2 = \alpha + t_2$.
 
Finally for (iv), let $V(G_2) =  \lbrace  x_1,x_2,\cdots,x_{s_1}\rbrace \cup \lbrace y_1,y_2,\cdots,y_{t_1}\rbrace \cup \lbrace z_1 \rbrace \cup \lbrace u_1 \rbrace$.  Then the resolving set of $G_2$ contains the resolving set of $G_1$ along with one of the vertex from $\{z_1, u_1\}$ and therefore, $\beta(G_2) =\alpha + 1$.
 \end{prf}
 \label{x}
\end{prop}
\begin{thm} 
For each $s_i>1$, $t_i>1$ with $1\leq i\leq k$,  the metric dimension of a graph $G$ realized by a code of the type $0^{s_1}1^{t_1}0^{s_2}1^{t_2}\cdots0^{s_k}1^{t_k}$ is $\sum\limits_{i=1}^k (s_i~+~t_i~-~2)$.
\begin{prf}
Apply Proposition \eqref{x} to each string of the form $(0^{s_3}1^{t_3})$, $(0^{s_4}1^{t_4})$, $\cdots$, $(0^{s_k}1^{t_k})$. Then the resolving set of $G$ contains exactly $s_i~-1+~t_i~-~1$ number of vertices. Thus, $\beta(G)~=~ \sum\limits_{i=1}^k (s_i~+~t_i~-~2)$.
\end{prf}
\label{y} 
\end{thm}
Next, we provide bounds for the metric dimensions of a graph $G$ realized by any code of the type $0^{s_1}1^{t_1}0^{s_2}1^{t_2}\cdots0^{s_k}1^{t_k}$. We also present families of graphs for which the bounds are attained.
\begin{thm}
Let $G$ be a graph realized by a code $0^{s_1}1^{t_1}0^{s_2}1^{t_2}\cdots0^{s_k}1^{t_k}$. Then the following hold, $$\sum\limits_{i=1}^k (s_i + t_i - 2)\leq \beta(G)\leq \sum\limits_{i=1}^k (s_i+t_i-1).$$ If  $s_i > 1$ and $t_i>1$ with $1\leq i\leq k$, then the lower bound is attained and the upper bound is attained if and only if $s_i = 1$ and $t_i\geq 1$.
\begin{prf}
By Lemma \eqref{1}, the resolving set  of every graph realized by a string $(0^{s_i}1^{t_i})$ contains atleast $(s_i+t_i-2)$ vertices. Therefore, $\beta(G)\geq \sum\limits_{i=1}^k (s_i~+~t_i~-~2)$ and the equality follows from Theorem \eqref{y}.\\
$~~~~~$By Lemma \eqref{1}, we see that for each $i$, $1\leq i\leq k$ there are at most two vertices in each string $(0^{s_i}1^{t_i})$ (one from each partite set) which does not belong to $\mathcal{W}$ of $G$. If $s_i = 1$ for all $i$, then by Theorem \eqref{y}, the set $\mathcal{W}$ does not contain any vertex from the independent part of the string $(0^{s_i}1^{t_i})$. This implies that $\mathcal{W}$ contains all vertices of the complete part. However, if $s_i\geq 2$, then $\mathcal{W}$ contains atleast one vertex from the independent part of the string $(0^{s_i}1^{t_i})$, which resolves any pair of vertices from two consecutive complete parts say $J_i$ and $J_{i+1}$ of $G$. Note that if $s_i = 1$ for each $i$, then $|\mathcal{W}|$ is maximum (that is, when every string in the graph $G$ is of the form $01^{t_i}$, where $t_i \geq 2$), that is, $|\mathcal{W}|\leq \sum\limits_{i=1}^k(s_i+t_i-2+1) =\sum\limits_{i=1}^k(s_i+t_i-1)$. Moreover, if for each $i$, $s_i=1$ and $t_i\geq 1$, then the equality follows  directly from Theorem \eqref{y} and consequently, $\beta(G)= \sum\limits_{i=1}^k (1+t_i-1)= \sum\limits_{i=1}^kt_i$.\\
 $~~~~~~~$Conversely, suppose that there is atleast one $s_i > 1$. Without loss of generality say it is $s_1$. Then the graph $G$ realizes a code $0^{s_1}1^{t_1}01^{t_2}\cdots01^{t_k}$. By Lemma \eqref{1}, we have seen that if the cardinality of the independent part  in any string $0^{s_i}1^{t_i}$ is greater than $1$, then the resolving set contains exactly $s_i+t_i-2$ vertices. Using this argument, the metric dimension of $G$ is 
 $(s_1+t_1-2)+(s_2+t_2-1)+\cdots+(s_k+t_k-1)=\sum\limits_{i=1}^k (s_i+t_i-1)-1<\beta(G)$, a contradiction.
\end{prf} 
\label{i}
\end{thm}
It is pertinent to mention that from preceding results, the metric dimension $\beta(G)$ of a  graph $G$ realized by any code of the type $0^{s_1}1^{t_1}0^{s_2}1^{t_2}\cdots0^{s_k}1^{t_k}$ is dependent on the powers of bits $0$ and $1$ of strings $(0^{s_i}1^{t_i})$, that is, by varying the powers of bits in the code, $\beta(G)$ also varies. In fact $\beta(G)$ gets effected  if and only if $s_i=1$ or $t_i=1$. As seen above in Theorem \eqref{i}, the upper bound for $\beta(G)$ is attained if and only if $s_i = 1$ and $t_i\geq 1$ for each $i$ and the lower bound is attained if for each $i$, $s_i > 1$ and $t _i > 1$. Below, we present an algorithm to determine the metric dimension of any graph $G$ realized by a code of the type $0^{s_1}1^{t_1}0^{s_2}1^{t_2}\cdots0^{s_k}1^{t_k}$. The input is a binary code of length $k\geq 3$, and in return we obtain $\beta(G)$ .

\textbf{Algorithm-(a) (method to compute the metric dimension of threshold graphs)}

\noindent\textbf{INPUT:} $(0^{s_1}1^{t_1})(0^{s_2}1^{t_2})\cdots(0^{s_k}1^{t_k}).$\\
\textbf{OUTPUT:} $\beta(G)$.\\
\textbf{Step-1.} If $i=2$, $1\leq i \leq k$, then $\beta(G)=\alpha_1$ (computed in Proposition \eqref{x}). For $i=3$, go to Step-2.\\
\textbf{Step-2.} If $s_3>1$ and $t_3\geq 1$, then $\beta(G)= \alpha_1+(s_3+t_3-1)=\alpha_2$, only if $s_2>1,~t_2=1$ or $s_2>1$ and $t_1=s_1=t_2=1$, otherwise, apply   Proposition \eqref{x}. For $i=4$, go to Step-3.\\
\textbf{Step-3.}  If $s_4>1$ and $t_4\geq 1$, then $\beta(G)= \alpha_2+(s_3+t_3-1)=\alpha_3$, only if $s_3>1,~t_3=1$ or $s_2>1$ and $t_2=s_3=t_3=1$ or $s_1>1$ and $t_1=t_2=t_3=s_2=s_3=1$, otherwise, apply Proposition \eqref{x}. For $i=5$, go to Step-4.\\ \textbf{Step-4.} Follow the same iteration for $i = 5, 6, \cdots, k-1$ and suppose that for $i = k-1$, $\beta(G)=\alpha_{k-2}$. Finally, for $i=k$\, go to Step-5.\\
\textbf{Step-5.} If $s_k>1$ and $t_k\geq 1$, then $\beta(G)= \alpha_{k-2}+(s_k+t_k-1)=\alpha_{k-1}$, only if $s_{k-1}>1$ and $t_{k-1}=1$ or for any $i$, $1\leq i\leq k-2$, $s_i>1$, $t_i=1$ and $s_{i+1}=s_{i+2}=\cdots=s_{k-1}=t_{i+1}=t_{i+2}= \cdots =t_{k-1}=1$, otherwise, apply Proposition \eqref{x}.\\
\textbf{Step-6.} END.
 
It is not difficult to see that the algorithm-(a) is linear. Indeed for each $i$, $1\leq i\leq k$, $\beta(G)$ is computed on the basis of previous iteration and we have at most $n$ iterations such that each one is performed with $O(1)$ operations, which gives $O(n)$ for these steps.

We notice that by removing all the edges in the complete part of a threshold graph, the resultant graph is a chain graph. As discussed above, chain graphs can also be generated by a binary code of the type $0^{s_1}1^{t_1}0^{s_2}1^{t_2}\cdots0^{s_k}1^{t_k}$. The metric dimension of a threshold graph $0^s1^t$ is either $s+t-2$ or $s+t-1$. However, the metric dimension of its corresponding chain graph of order greater or equal to three is always $s+t-2$.

In the following result, we determine the metric dimension of a chain graph realized by a code $0^{s_1}1^{t_1}0^{s_2}1^{t_2}$.

\begin{lem}
Let $H_1$ be a chain graph realized by $0^{s_1}1^{t_1}$ with $\beta(H_1) = \alpha'$ and let $H_2$ be a chain graph realized by $0^{s_1}1^{t_1}0^{s_2}1^{t_2}$. Then the following hold.\\
(i) For each $s_2 > 1$ and $t_2 > 1$, 
\begin{equation*}
\beta(H_2)=\begin{cases}
         \alpha'+(s_2+t_2-1)  & if ~ s_1 > 1 ~\mbox{and}~ t_1=1 \\
         
         \alpha'+(s_2+t_2-2) & otherwise.
     \end{cases}
\end{equation*}
(ii) For $s_2>1$ and $t_2=1$,
\begin{eqnarray*}
\beta(H_2)=\begin{cases}
         \alpha'+s_2  & if ~ s_1 > 1 ~\mbox{and} ~t_1 = 1 \\
         
         \alpha'+(s_2-1) & otherwise.
     \end{cases}
\end{eqnarray*}
(iii) For $s_2 = 1$ and $t_2 > 1$, 
$\beta(H_2)=\alpha' + t_2$.\\\
(iv) For each $s_2 = t_2 = 1$,
\begin{equation*}
\beta(H_2)=\begin{cases}
         1  & if ~ s_1=t_1=1 \\
         
         \alpha'+1 & otherwise.
     \end{cases}
\end{equation*}

\begin{prf}
We partition the vertex set of $H_2$ as, $V(H_2)=U\cup V=\lbrace{u_1,u_2,\cdots,u_{s_1}\rbrace}\cup\lbrace{v_1,v_2,\cdots,v_{t_1}\rbrace}\cup\lbrace{u'_1,u'_2,\cdots,u'_{s_2}\rbrace}\cup\lbrace{v'_1,v'_2,\cdots,v'_{t_2}\rbrace}=U_1\cup V_1\cup U_2\cup V_2$, where $U_1, U_2$ represents classes of white vertices and $V_1, V_2$ represents classes of black vertices in $H_2$. Then for every $i, j, k$ and $l$,  $d(u_i, u'_k) = d(v_j, v'_l) = 2$,  $d(u_i, v_j)= d(u_i, v'_l)=d(u'_k,v'_l) = 1$ and $d(v_j, u'_k) = 3$.

To prove (i), suppose $s_1 > 1$ and $t_1 = 1$. Then, $V(H_2)=\lbrace{u_1,\cdots,u_{s_1}\rbrace}\cup\lbrace{v_1\rbrace}\cup\lbrace{u'_1,\cdots,u'_{s_2}\rbrace}\cup\lbrace{v'_1,v'_2,\cdots,v'_{t_2}\rbrace}$. From the structure of the graph $H_2$, we see that there are vertices $x \in U_1\cup U_2$, $y\in V_2$ and $z\in V_1$ such that $d(x|\mathcal{W}) \neq d(y|\mathcal{W}) \neq d(z|\mathcal{W})$, where $\mathcal{W} = V(H_2)\setminus \{x, y, z\}$. Therefore, $\beta(H_2)= s_1-1+s_2+t_2-1=\alpha'+(s_2+t_2-1)$,  where $\alpha'=\beta(H_1)$.\\
 $~~~~~$If $s_1=1$ and $t_1>1$ , then $V(H)= \lbrace{u_1\rbrace}\cup\lbrace{v_1,\cdots,v_{t_1}\rbrace}\cup\lbrace{u'_1,\cdots,u'_{s_2}\rbrace}\cup\lbrace{v'_1,\cdots,v'_{t_2}\rbrace}$. Using the same argument as above, the resolving set of $H_2$ is given as, $\mathcal{W}=\lbrace{v_1,\cdots, v_{t_1-1}, u'_1,\cdots,u'_{s_2-1},v'_1,\cdots,v'_{t_2-1}\rbrace}$. Thus, $\beta(H)=t_1-1+s_2-1+t_2-1=\alpha'+(s_2+t_2-2)$.
 
 For (ii), it is given that $s_2 > 1$ and $t_2 = 1$. Then for $s_1 > 1$ and $t_1 = 1$, $V(H_2)=\lbrace{u_1, u_2,\cdots, u_{s_1}\rbrace}\cup\lbrace{v_1\rbrace}\cup\lbrace{u'_1, u'_2,\cdots,u'_{s_2}\rbrace}\cup\lbrace{v'_{1}\rbrace}$. There are vertices $x \in U_1\cup U_2$, $y\in V_2$ and $z\in V_1$ such that $d(x|\mathcal{W}) \neq d(y|\mathcal{W}) \neq d(z|\mathcal{W})$, where $\mathcal{W} = V(H_2)\setminus \{x, y, z\}$. Therefore, $\beta(H)= |U_1| + |U_2| - 1= s_1 - 1 + s_2=\alpha'+s_2$.\\
 $~~~~~$If $s_1 = 1$ and $t_1 > 1$, then using the same argument as above, the resolving set of $H_2$ is, $\mathcal{W}=\lbrace{v_1,v_2,\cdots,v_{t_1-1},u'_1,u'_2,\cdots,u'_{s_2-1}\rbrace}$ and consequently, $\beta(H_2) = t_1-1+s_2-1=\alpha'+s_2-1$.
 
To prove (iii), let $s_1>1$ and $t_1=1$. Then, $V(H_2)=\lbrace{u_1, u_2,\cdots, u_{s_1}\rbrace}\cup\lbrace{v_1\rbrace}\cup\lbrace{u'_1\rbrace}\cup\lbrace{v'_1, v'_2,\cdots, v'_{t_2}\rbrace}$. The proof follows from (ii).

For (iv), it is given that $s_2 = t_2 = 1$, suppose that $s_1= t_1 = 1$. Then the graph $H_2$ is a path on four vertices and therefore, $\beta(H_2)=1$. If $s_1 > 1$ and $t_1 = 1$, then the basis of $H_2$ contains all vertices of $H_2$ except one vertex of $U_1$. Thus, $\beta(H_2)= s_1-1 + 1=\alpha' + 1$.\\ 
\end{prf}
\label{b} 
  \end{lem}
\begin{rem}
By Lemma \eqref{b}, the metric dimension of any chain graph $H$ realized by $(0^{s_1}1^{t_1})(0^{s_2}1^{t_2})\\
\cdots(0^{s_i}1^{t_i})$ can be obtained iteratively for each $i\geq 3$. In fact, for $3\leq i\leq k$, the  Algorithm-(a) can be applied to determine the metric dimension of $H$, note that for $i = 2$, $\beta(H)=\alpha_2$ (Lemma \eqref{b}). The Algorithm-(a) establishes a relationship between $\beta(H)$ and $\beta(G)$ and the equality $\beta(H) = \beta(G) -1$ or $\beta(H) = \beta(G)$ follows, where $G$ is threshold graph.
\end{rem}

Now, we provide an example which illustrates the above theory and Algorithm-(a) related to $\beta(H)$ and $\beta(G)$.
 
\begin{exm}
 Let $G$ be a threshold graph and $H$ be a chain graph generated by a code of the type $01^20^310^21^20^31^4010^31^4$. By Algorithm-(a), we see that $\beta(G) = 18$ whereas $\beta(H) = 17$. Moreover, for $s_1\geq 2$, if  $G$ and $H$ are realized by $0^{s_{1}}10^{s_2}10^{s_3}10^{s_4}10^{s_5}10^{s_6}10^{s_7}10^{s_8}1$, then $\beta(H)= \beta(G) = \big(\sum\limits_{i = 1}^{8}s_i\big) - 1$. 
 \end{exm}
 
 \section{Threshold dimension, restricted threshold dimension of a threshold graph}
 
 The metric basis determines the location of landmarks in a graph with a fixed number of vertices. Chartrand et. al \cite{GLMO} and Khuller et. al \cite{SBA2} previously posed a question that how the metric dimension of a graph is related to its subgraphs? However, in \cite{LMO}, the authors have asked that how much the metric dimension of a graph can be reduced by adding edges? Equivalently, for a given graph $G$, what is the smallest metric dimension among all graphs that contain $G$ as a spanning subgraph? Recall that the \textit{threshold dimension} of a graph $G$, denoted by $\tau(G)$, is defined as,
 \begin{center}$\tau(G)=min\lbrace \beta(H): \mbox{H contaning G as a spanning subgraph}\rbrace$.\end{center}
 A graph $G$ is called \textit{irreducible} if $\beta(G) = \tau(G)$, otherwise it is \textit{reducible}. In this section, a graph of order less or equal to $4$ is always irreducible.
 
 \begin{lem}
For $s\geq 3$, let $G$ be a graph realized by $0^s1^t$. If $s \leq 2^k- 1 + k,$ then $\tau(G) = t - 1 + k$, where $k\geq 2$ is the minimum number of isolated vertices taken from the independent part of $G$. 
 
\begin{prf}
Let $V(G)=I\cup J$, where $I$ and $J$ represents the independent and complete parts of $G$ with cardinalities $s$ and $t$ respectively. The only way to add edges in $G$ is to add edges between the isolated vertices in $I$. Choose $W\subset I$, a proper subset of $I$ with cardinality $k\geq 2$. Then it can be verified that there are $2^{k} - 1$ vertices of $I$ such that their distance vectors are distinct with respect to $W$. In fact, the distance vectors are of the form $(1, 1, \cdots, 1, 2),\cdots, (1, 2, 2, \cdots, 2, 1),\cdots, (2, 1,\cdots, 1), \cdots, (2, 2, \cdots, 2)$. Therefore the resolving set of $G$ contains $|J|-1$ and $|W|$ vertices. Thus, $\tau(G) = t - 1 + k$.  
\end{prf}
\label{cc}
\end{lem}
\begin{exm}
Let $G$ be a graph realized by $0^{2359}1^{15}$. Then for $k=11$, $\tau(G)=15-1+11 = 25$ and $2359 - 11 \leq  2^{11}-1$. However, if $k = 10$, then it can be verified that for $G$, Lemma \eqref{cc} does not hold.
\end{exm}

\begin{thm} \label{23}
Let $G$ be a threshold graph realized by $0^{s_1}1^{t_1}0^{s_2}1^{t_2}$. Then  one of the following holds.\\
(i) If $s_1+t_1+s_2-k_1\leq 2^{k_1}-1$, then $\tau(G) = (t_2-1)+k_1$. \\
(ii) If $t_1\geq 2^{s_2}-1$ and $s_1-k_2\leq 2^{s_2+k_2}-2^{s_2}$, then $\tau(G) = (t_2-1)+s_2+m_1+k_2$. \\
(iii) If $t_1< 2^{s_2}-1$ and $s_1-m_2-k_2\leq 2^{s_2+k_2}-2^{s_2}$, then $\tau(G) = (t_2-1)+s_2+k_2$,\\ 
$~~~~~~~\mbox{where}, m_1=t_1-2^{s_2}+1$, $m_2=2^{s_2}-1-t_1$ and $k_1$, $k_2$ are positive integers.
\end{thm}
\begin{prf}
Let $V(G)=U\cup V$, where $U=U_1\cup U_2$ and $V=V_1\cup V_2$ represents the independent and complete parts of $G$ with cardinalities $s_1+s_2$ and $t_1+t_2$. Clearly, from the structure of $G$, the resolving set of $G$ contains atleast $|V_2|-1$ vertices.\\
To prove (i), choose $W\subseteq U_2$, a proper subset with cardinality $k_1\geq 2$. Then by the similar argument as in Lemma \eqref{cc}, we have the required result.

If $|V_1|\geq 2^{|U_2|}-1$, then $ t_1-2^{s_2}+1\geq 0$. Note that $m_1 = t_1-2^{s_2}+1$ vertices from $V_1$ have same distance vectors with respect to $U_2$ after adding $s_2(2^{s_2-1}-1)$ edges between $V_1$ and $U_2$. Further adding of edges between vertices of $V_1$ is not possible, since $V_1$ is complete. It remains to add edges between the vertices of partite sets $U_1$ and $U_2$. As before, choose $W\subseteq U_1$, a proper subset with cardinality $k_2$. Then the resolving set of $G$ contains exactly $|V_2|-1+|U_2|+|V_1| - 2^{|U_2|}+1 +k_2$ vertices, provided $|U_1|-k_2\leq 2^{|U_2|+k_2}-2^{|U_2|}$. Thus, $\tau(G)=(t_2-1)+s_2+m_1+k_2$, whenever $s_1-k_2\leq 2^{s_2+k_2}-2^{s_2}$. Therefore, (ii) follows.
  
If $|V_1|<2^{|U_2|}-1$, then $m _2 = t_1-2^{s_2}+1<0$ and there is nothing to prove. Suppose $|U_1|\leq |m_2|$, then it gets reduced to (i) proven above. If  $|U_1|> |m_2|$, then as in (ii), we take a smallest subset of $U_1$ with cardinality as $k_2$ satisfying $m_2-k_2\leq 2^{s_2+k_2}-2^{s_2}$. Therefore resolving set of $G$ contains $(t_2)-1+s_2+k_2$ vertices, whenever $m_2-k_2\leq 2^{s_2+k_2}-2^{s_2}$. Thus, (iii) follows. 

\end{prf}

\begin{thm} \label{24}
Let $G$ be a threshold graph realized by $0^{s_1}1^{t_1}0^{s_2}1^{t_2}0^{s_3}1^{t_3}$. Then one of the following holds.\\
(i) If $\sum\limits_{i=1}^3 (s_i~+~t_{3-i})-k_1\leq 2^{k_1}-1$, then $\tau(G)=(t_3-1)+k_1$. \\
(ii) If $t_2\geq 2^{s_3}-1$ and $s_1+t_1+s_2-k_2\leq 2^{s_3+k_2}-2^{s_3}$, then $\tau(G)=(t_3-1)+s_3+\big(t_2-(2^{s_3}-1)\big)+k_2$.\\
(iii) If $t_2< 2^{s_3}-1$, $\big(t_1-(2^{s_3}-1-t_2)\big)\geq (2^{s_2+s_3}-2^{s_3})$ and $s_1+s_2-k_3\leq 2^{s_3+k_3}-2^{s_3}$, then $~~~~~~~\tau(G)=(t_3-1)+s_3+(t_1+t_2+1-2^{s_2+s_3})+k_3$. \\
(iv) If $t_1+t_2<2^{s_3}-1$, $\sum\limits_{i=1}^2 (s_i~+~t_i)\geq 2^{s_3}-1$ and  $\sum\limits_{i=1}^2 (s_i~+~t_i)-k_3\leq 2^{s_3+k_3}-2^{s_3}$, then $~~~~~~~\tau(G)=(t_3-1)+s_3+k_3$.\\
The numbers $k_3, k_2$ and $k_1$ are any positive integers.
\end{thm}
\begin{prf}
The inequalities of the above type are extensions of the inequalities given in Theorem  \eqref{23}, therefore the result can be proven by the similar approach.

\end{prf}
\begin{rem} By Theorems \{\eqref{23}, \eqref{24}\}, we can determine the threshold dimension of any threshold graph $G$ realized by a code of the type $0^{s_1}1^{t_1}0^{s_2}1^{t_2}\cdots0^{s_k}1^{t_k}$. Suppose for each $i$, $s_i\geq t_i$ and $\sum\limits_{i=1}^k (s_i~+~t_{k-i})-k_1\leq 2^{k_1}-1$, where $1\leq i\leq k$ and $k_1$ a positive integer. Then $\tau(G)=(t_k-1)+k_1$. Moreover, if $\sum\limits_{i=1}^{k -1} < 2^{s_k}-1$, $\sum\limits_{i=1}^{k-1} (s_i~+~t_i)\geq 2^{s_k}-1$ and  $\sum\limits_{i=1}^{k-1} (s_i~+~t_i)-r\leq 2^{s_k+r}-2^{s_k}$, then $\tau(G)=(t_k-1)+s_k+r$, $r$ is some positive integer. Therefore, based on the conditions of the powers of bits $0$ and $1$ of the code we can establish the threshold dimension of any threshold graph.\\
\end{rem}

\subsection{Restricted threshold dimension}

As defined in the preceding section, the threshold dimension of a threshold graph $G$ is the minimum metric dimension among all graphs that contain $G$ as a spanning subgraph. However, the resulting graph with the least metric dimension may or may not be threshold. Related to the threshold dimension of $G$, there is another graph theoretical invariant called as the \textit{restricted threshold dimension} of $G$, denoted by $\tau_r(G)$, which is defined as the smallest metric dimension among all threshold graphs which contains $G$ as a spanning subgraph, that is, $$\tau_r(G)=min\lbrace \beta(H): \mbox{H is a threshold graph contaning G as a spanning subgraph}\rbrace.$$
This section is devoted to study the restricted threshold dimension of the threshold graph.

\begin{prop}
Let $G$ be a threshold graph realized by $0^s1^t$.  Then the following hold, $$\tau_r(G) = \lceil\frac{s}{2}\rceil+(t-1),$$
where $s\geq 4$.
\begin{prf}
We partition the vertex of $G$ as, $V(G)~ =~ U\cup V$. The set $U~ = \lbrace u_1,u_2,\cdots,u_s\rbrace$,   represents the independent part of $G$ and $V~ = \lbrace v_1,v_2,\cdots,v_t\rbrace$, represents the complete part of $G$. By Lemma \eqref{1}, $\beta(G)=s+t-2$. The only way to add edges in $G$ such that $G$ continues to remain a threshold is to add edges between the vertices in $U$. Suppose $u_1 \sim \lbrace u_2,u_3,\cdots,u_s\rbrace$, that is $u_1$ is adjacent to each of the vertex from $\lbrace u_2,u_3,\cdots,u_s\rbrace$. Then the resulting graph $H$ is the threshold graph realized by$0^{s-1}1^{t+1}$ with $\beta(G)=\beta(H)$. Now, suppose $u_1\sim u_2$ and $\{u_1, u_2\} \sim \lbrace u_3,u_4,\cdots,u_s\rbrace$. Then the resulting threshold graph is some $H$ realized by $0^{s-2}1^{t+2}$ with $\beta(G)=\beta(H)$. We continue the same process of adding edges till the power of a bit $0$ is reduced to two and the resulting graph is a threshold graph realized by $0^21^{s+t-2}$.  The total number of such graphs realized by strings of length $1$ are ${{s}\choose {1}}-2$. Further, let $H$ be a threshold graph realized by $0^{s_1}1^{t_1}0^{s_2}1^{t}$ containing $G$ as a spanning subgraph. Note that $H$ is obtained from $G$ by making $t_1$ dominating vertices in $U$ and partition $s$ isolated vertices of $U$ into $s_1$ and $s_2$ isolated vertices. The total number of such graphs is ${{s}\choose {3}}$. Using the same procedure, it can be verified that if $s$ is even, then there are exactly ${{s}\choose {1}}-2+{{s}\choose {3}}+{{s}\choose {5}}+\cdots+{{s}\choose {s-1}}$ different threshold graphs with strings of lengths $1,2,\cdots, \frac{s}{2}$. If $s$ is odd, there are  ${{s}\choose {1}}-2+{{s}\choose {3}}+{{s}\choose {5}}+\cdots+{{s}\choose {s-2}}+1$ different threshold graphs with lengths $1,2,\cdots, \lceil\frac{s}{2}\rceil$.\\
By Theorem \eqref{i}, among all $2^{s-1}-2$ threshold graphs with strings of lengths less or equal to $\lfloor\frac{s}{4}\rfloor$, the metric dimension is greater than $\lceil\frac{s}{2}\rceil$. Using the Algorithm-(a), the metric dimension of all these $2^{s-1}-2$ threshold graphs with $G$ as a spanning subgraph lies in $[\lceil\frac{s}{2}\rceil, s+t-2]$, where the maximum is attained for graphs with strings of length $1$ and the minimum is attained if the graph with string is of length atleast $\lfloor\frac{s}{4}\rfloor+1$.  Furthermore, the graphs obtained by adding the maximum number of edges with smaller metric dimension are clearly dependent on the cardinality of $U$ alone. So, we consider the following cases.\\
{\bf Case - I.} $s\in \{4k : k\geq 1\}\cup \{4k +2 : k\geq 1\}$, then following the above computations, there is a unique graph $H$ realized by $0^210101\cdots0101^t$ with strings of length $\frac{s}{2}$. In fact $H$  is obtained by adding $2k(2k-1)$ or $2k(2k+1)$ edges with metric dimension $\underbrace{1+1+\cdots+1}_{\frac{s}{2}-1}+t=\frac{s}{2}+t-1$\\
{\bf Case - II.}  $s\in \{4k+1 : k\geq 1\} \cup\{4k-1 : k\geq 2\}$, then there is a unique graph $H=0^210101\cdots0101^{t+1}$  with strings of length $\frac{s-1}{2}$. The graph $H$ is obtained by adding $2k(2k+1)$ or $2k(2k-1)$ edges with metric dimension $\underbrace{1+1+\cdots+1}_{\frac{s-1}{2}-1}+(t+1)=\frac{s-1}{2}+t=\lceil\frac{s}{2}\rceil+t-1$.\\
Thus, $\tau_r(G)=\beta(H)=\lceil\frac{s}{2}\rceil+t-1$.
\end{prf}
\label{rt}
\end{prop}

\begin{thm}
Let $G$ a threshold graph realized by $0^{s_1}1^{t_1}0^{s_2}1^{t_2}\cdots0^{s_k}1^{t_k}$. Then for each $i$, $ \tau_{r}(G) = \sum\limits_{i=1}^{k}\big(\lceil\frac{s_i}{2}\rceil+(t_i-1)\big)$, where $1 \leq i \leq k$, $s_i\geq 4$ and $t_i\geq 2$. 
\end{thm}
\begin{prf}
By Proposition \eqref{rt}, the string $0^{s_1}1^{t_1}$ can be expressed in two ways only. If $s_1$ is even, then the string can be uniquely written as $0^210101\cdots0101^{t_1}$. On the other hand, if $s_1$ is odd, then it can be written as $0^210101\cdots0101^{t_{1}+1}$. Using this approach to each string $0^{s_i}1^{t_i}$, $2\leq i\leq k$ of the graph $G$. Since $G$ has $k$ strings, therefore, the metric dimension of $H$ is $\big(\lceil\frac{s_1}{2}\rceil+(t_1-1)\big)+\big(\lceil\frac{s_2}{2}\rceil+(t_2-1)\big)+\cdots+\big(\lceil\frac{s_n}{2}\rceil+(t_n-1)\big)=\sum\limits_{i=1}^{n}\big(\lceil\frac{s_i}{2}\rceil+(t_i-1)\big)=\tau_r(G)$.
\end{prf}
\begin{rem}
For any threshold graph $G$ realized by  $0^{s_1}1^{t_1}0^{s_2}1^{t_2}\cdots0^{s_k}1^{t_k}$, $\beta(G)=\tau_r(G)$, where $1\leq i \leq k$, $s_i\leq 3$ and $t_i\geq 2$. If $s_i \geq 4$ and $t_i = 1$, then $\tau_r(G)=\lceil\frac{s_1}{2}\rceil+(\lceil\frac{s_2}{2}\rceil+1)+\lceil\frac{s_3}{2}\rceil+(\lceil\frac{s_4}{2}\rceil+1)+\cdots+\big((\lceil\frac{s_k}{2}\rceil+1)$ or $\lceil\frac{s_k}{2}\rceil\big)$ (depending on the length of the binary generating code). If $s_i\leq 3$ and $t_i=1$, then $\tau_r(G)=\sum\limits_{i=1}^{k}\lceil\frac{s_i}{2}\rceil$. Summing up all the theory related to the restricted threshold dimension of a threshold graph $G$ realized by $0^{s_1}1^{t_1}0^{s_2}1^{t_2}\cdots0^{s_k}1^{t_k}$, $\tau_r(G)$ can be obtained by the following procedure.\\
(i) If $s_i \in \{4k : k\geq 1\}\cup \{4k +2 : k\geq 1\}$, then write each of the strings $0^{s_i}1^{t_i}$ of $G$ as $0^210101\cdots0101^{t_i}$. If $s_i \in\{4k+1 : k\geq 1\} \cup\{4k-1 : k\geq 2\}$, then write $G$ as $0^210101\cdots0101^{t_{i}+1}$. In each case the strings of $G$ are of the form $(0^31^{t_i})$, we write it as  $0101^{t_i}$, and continue the process untill the power of a bit $0$ in the resulting graph is reduced to not more than $2$.\\
(ii) If the code of the resulting graph $H$ is free from the string $0^{s}1^{t}$, $t\geq 2$ and strings $(0^21)$ are repeated $n$ number of times, then replace $0^21\longrightarrow (01^2)$ at even places of the code. If in the code of $H$, the string $(0^21)$ is repeated $n$ number of times between strings $0^{s}1^{t}$, $t\geq 2$, then again replace $(0^21)\longrightarrow (01^2)$ at even places of the code. Finally, if the string $(0^21)$ is repeated $n$ number of times before or after the string $0^{s}1^{t}$, $t\geq 2$, then     follow the same process as mentioned. Therefore, by the Algorithm-(a), $\beta(H)=\tau_r(G)$.
\label{rs}
\end{rem}
The following example illustrates Proposition \eqref{rt} and Remark \eqref{rs}.
\begin{exm}
Let  $G$ be a threshold graph realized by $0^31^20^81^20^510^61^40^71^2$. Then by Algorithm-(a), $\beta(G) = 29$. Further, by applying Proposition\eqref{rt} to each of the strings of $G$, the code of the resulting threshold graph $H$ is $(0101^2)(0^21010101^2)(0^2101^2)(0^210101^4)(0^210101^3)$ and $H$ contains $G$ as a spanning subgraph. Therefore, again by the Algorithm-(a), $\tau_{r}(G)= 22 = \beta(H)$.
 However, if $G=0^210^410^610^81$, then $\beta(G)=19$ and the code of the resulting graph $H$ is $(0^21)(0^2101)(0^210101)(0^21010101)$. Note that the code of $H$ contains the string $0^21$ four times. Therefore, by Remark \eqref{rs} replace $0^21\longrightarrow 01^2$ at even places. Consequently, $\tau_{r}(G)= 12 = \beta(H)$.
\end{exm}
\section{$\lambda$-chromatic number}

Let $j$ and $k$ be two positive integers, an $L(j, k)$-coloring of a graph $G$
is a function from $V(G)$ to the set of non-negative integers such that $|f(x) - f(y)|\geq j$, whenever $x, y\in V(G)$ are adjacent and $|f(x) - f(y)| \geq k$, whenever $x, y\in V(G)$ are at distance two apart. The difference between the
maximum and minimum values of $f$ is called its span denoted by $span (f)$. If $f$ is a $L(j, k)$-coloring of $G$ with a minimum value $\delta$ (say), then the function $\mathfrak{g}$ defined on $V(G)$ by $\mathfrak{g}(v) = f(v)$ is also an $L(j, k)$-coloring with a minimum value of $0$ and a maximum value of $span(f) = span(\mathfrak{g})$. An \textit{$L(j, k)$-chromatic number}, denoted by $L_{j,k}(G)$, is the minimum span over all $L(j, k)$-colorings of $G$, and the minimal coloring is the $L(j, k)$ coloring that corresponds to it. For the cases $j = 2$ and $k = 1$, an $L(2, 1)$-coloring problem was first studied by Griggs and Yeh $\cite{Gy}$. An $L(2, 1)$- coloring that uses colors from the set $\lbrace 1, 2, \cdots, s\rbrace$ (not necessarily all) is called an $s$-coloring. The minimum $s$, so that $G$ has a $s$-coloring is called the \textit{$\lambda$- chromatic number} of $G$ and is denoted by $\lambda(G)$. If the color $h$, $0 < h < s$, is not used in an $s$-coloring of $G$, then $h$ is called the \textit{hole} of the coloring. The $L(j, k)$-coloring problem and $L(2, 1)$-coloring problem have been extensively studied in $\cite{Gm,KR}$.

In this section, we study the $\lambda$-chromatic number of threshold and chain graphs. In the following result, we exactly determine the $\lambda$- chromatic number of a graph with diameter two.

\begin{lem}
Let $G$ be a diameter two graph with   $\lambda(G) = s$. Then the following hold, $$\lambda(\mathcal{G}_{m,n})=s+\frac{|m-k|+(m-k)}{2}+2n,$$ where $\mathcal{G}_{m,n}$ is the graph obtained from $G$ by adding $m$ isolated and $n$ dominating  vertices, and $k$ is the number of holes in the $s$-coloring of $G$.

\begin{prf}
Let $f$ be the minimum $L(2, 1)$-coloring of $G$ with  $\lambda(G)=s$ and let $\lbrace h_1, h_2,\cdots, h_k\rbrace$ be its holes. We obtain $\mathcal{G}_{m, n}$ from $G$ by adding $\lbrace x_1, x_2, \cdots, x_m\rbrace$ isolated vertices and $\lbrace y_1, y_2, \cdots, y_n\rbrace$ dominating vertices. Define the $L(2, 1)$-coloring $g$ on $\mathcal{G}_{m,n}$ as follows,
\begin{eqnarray*}
g(y_{i})& =& 2i-2, ~~~~1\leq i\leq n,\\
g(v) &=& f(v)+2n,~~v\in V(G),\\
g(x_i) &=& h_i+2n, ~~~~1\leq i\leq m. 
\end{eqnarray*}
If $m>k$ and $1\leq j\leq m-k$, then we define, $g(x_{k+j})=s+j+2n$.
$$\mbox{Therefore}, \max\limits_{v\in V(G_{m,n})} {g(v)}=s+\frac{|m-k|+(m-k)}{2}+2n.$$
Since the diameter of $G$ is two, therefore $f$ must be one-one and this implies $g$ is well defined $L(2, 1)$-coloring of $\mathcal{G}_{m,n}$. Further, $f$ is minimum, therefore, $\lambda(\mathcal{G}_{m,n})<s+\frac{|m-k|+(m-k)}{2}+2n$ is not possible and we conclude $\lambda(\mathcal{G}_{m,n})=s+\frac{|m-k|+(m-k)}{2}+2n$.
\end{prf}
\label{ln}
\end{lem}

The following result is an immediate consequence of Lemma \eqref{ln}.
\begin{thm}
Let $G$ be any threshold graph realized by $(0^{s_1}1^{t_1})(0^{s_2}1^{t_2})\cdots(0^{s_k}1^{t_k})$. Then $$\lambda(G)=2\sum\limits_{i=1}^{k}t_i+(s_1-1)+\frac{1}{2}{\sum\limits_{i=1}^{k-1}\big(|s_{i+1}-h_i|+(s_{i+1}-h_i)\big)}$$ where $h_i=t_i+\frac{1}{2}\big(|t_{i-1}-s_i|+(t_{i-1}-s_i)\big)$, $1\leq i\leq k$, denotes the number of holes in the minimum coloring of $G$.
\end{thm}

\begin{prf}
To prove the result, we use an iterative procedure on $k$. For $k=1$, let $\mathcal{G}$ be a graph realized by $(0^{s_1}1^{t_1})$ with $V(\mathcal{G})=X\cup Y=\lbrace x_1, x_2, \cdots, x_{s_1}\rbrace \cup \lbrace y_1, y_2, \cdots, y_{t_1}\rbrace$. Define the $L(2, 1)$-coloring $f$ on $V(\mathcal{G})$ as follows,
\begin{eqnarray*}
f(y_i) &=& 2i-2,~~~~~~~~ 1\leq i \leq t_1,\\
f(x_i)& =& 2t_1+(i-1),~~~ 1\leq i \leq s_1.
\end{eqnarray*}
Since $\mathcal{G}$ is threshold with diameter less than three, this implies that $f$ is one-one. Moreover, the clique number of $\mathcal{G}$ is $t_1 + 1$, therefore $f$ has $t_1$ number of holes. Thus, $f$ is minimal and $\lambda(\mathcal{G})=2t_1+s_1-1$.  Let $\mathcal{G}_{s_2,t_2}$ be a graph realized by $(0^{s_1}1^{t_1})(0^{s_2}1^{t_2})$, obtained from $\mathcal{G}$ by adding $s_2$ isolated and $t_2$ dominating vertices. Then by Lemma \eqref{ln}, $$\lambda(\mathcal{G}_{s_2,t_2})=2\sum\limits_{i=1}^{2}t_i+(s_1-1)+\frac{1}{2}\big(|s_2-h_1|+(s_2-h_1)\big),$$ where $h_1 = t_1$. We iterate further and finally, we obtain $\mathcal{G}_{s_k,t_k}= G$ realized by $(0^{s_1}1^{t_1})(0^{s_2}1^{t_2})\cdots(0^{s_k}1^{t_k})$. Again, by Lemma \eqref{ln}, $$\lambda(G)=2\sum\limits_{i=1}^{k}t_i+(s_1-1)+\frac{1}{2}{\sum\limits_{i=1}^{k-1}\big(|s_{i+1}-h_i|+(s_{i+1}-h_i)\big)}.$$ For $1\leq i\leq k$, the equality $h_i=t_i+\frac{1}{2}\big(|t_{i-1}-s_i|+(t_{i-1}-s_i)\big)$  follows directly from the definition of  minimal coloring $f$ of $G$.
\end{prf}

Recall that a chain graph satisfy the nesting property, that is, there exists a partition of $U = U_1\cup U_2\cup \cdots \cup U_l$ and $V = V_1\cup V_2\cup \cdots \cup V_l$ such that for each $i$, $1\leq i\leq l$, the neighbourhood of each vertex in $U_i$ is $V_1\cup V_2\cup \cdots \cup V_{l+1-i}$. If $|U_i|= m_i$ and $|V_i|= n_i$, where $i = 1, 2,\cdots, l$, then the another representation of $G$ is $G(m_1,\cdots, m_l ;  ~n_1,\cdots,n_l)$. Below, we present an algorithm for determining the $\lambda$-chromatic number of chain graphs.

\textbf{Algorithm-(b) (method for computing the $\lambda$-chromatic number of chain graph)}\\

\noindent\textbf{INPUT:} Chain graph $G(m_1,\cdots, m_l ;  ~n_1,\cdots,n_l)$\\
\textbf{Output:} $\lambda(G)$\\
{\bf Step-1.} For each vertices of the type $u_{i_{1}}\in U_1$, $u_{i_{2}}\in U_2, \cdots, u_{i_{l}}\in U_l$, where $1\leq i_1\leq m_1$, $1\leq i_1\leq m_2, \cdots,1\leq i_l\leq m_l$. Define $L(2, 1)$-coloring on $V(G)$ as,
\begin{eqnarray*}
f(u_{i_{1}}) &=& i_1-1, ~~~~~~~~~~ 1\leq i_1\leq m_1,\\ f(u_{i_{2}}) &=& m_1-1+i_2,  ~~~ 1\leq i_2\leq m_2,\\\vdots\\
f(u_{i_{l}}) &=& \sum\limits_{i=1}^{l-1}(m_i)-1+i_l~~~ 1\leq i_l\leq m_l.
\end{eqnarray*}
and for vertices of type $v_{i_{1}}\in V_1, v_{1_{2}}\in V_2,$ where $1\leq i_1\leq n_1$, define $L(2, 1)$-coloring as,
$$f(v_{i_{1}}) =\sum\limits_{i=1}^{l}m_i+i_1, ~~1\leq i_1\leq n_1 ~\mbox{and}~
f(v_{1_{2}}) =\sum\limits_{i=1}^{l}m_i.$$
The sum $\sum\limits_{i=1}^{l}m_i$ is the hole of $f$. To color the remaining vertices of the graph $G$, we go to Step-2 (note that by definition of the graph, coloring of vertices of $v_{i_{2}}\in V_2$, $2\leq i_2\leq n_2$, depends on the colors of vertices of $U_l$).\\
{\bf Step-2.} Compute $|V_2|-|U_l|$,  if  $s_1 = |V_2|-|U_l|\leq 0$, then go to Step-3.\\
{\bf Step-3.} If $r_1 = |V_2|-|U_l|> 0$, define, 
\begin{eqnarray*}
f({v_{{i_{2}}^{'}}})&=&\sum\limits_{i=1}^{l}(m_i)+n_1+{i_{2}^{'}}~~~~1\leq {i_{2}^{'}}\leq r_1. 
\end{eqnarray*}
{\bf Step-4.} Compute $|V_3|-|U_{l-1}|$,  if  $ s_2 = |V_3|-|U_{l-1}|\leq 0$, then go to Step-5.\\
{\bf Step-5.} If  $r_2 = |V_3|-|U_{l-1}|> 0$, then compute $r_2 + s_1$. If $s_3 = r_2+s_1\leq 0$, then go to Step-6.\\ 
{\bf Step-6.} If $r_3 = r_2 + s_1>0$,  define,
\begin{eqnarray*}
f({v_{{i_{3}}^{'}}})&=&\sum\limits_{i=1}^{l}(m_i)+n_1+{i_{3}^{'}}, ~~1\leq {i_{3}^{'}}\leq r_3,\end{eqnarray*} otherwise go to Step-7.\\ 
{\bf Step-7.} If no such $s_1$ exists, then define,
\begin{eqnarray*}
f({v_{{i_{3}}^{'}}})&=& \sum\limits_{i=1}^{l}(m_i)+n_1+r_1+{i_{3}^{'}}, ~~~1\leq {i_{3}^{'}}\leq r_3.\\
\end{eqnarray*}
{\bf Step-8.} Continue the same iterations.\\
{\bf Step-9.} End.

Next, we provide extreme values (bounds) for $\lambda$-chromatic number of chain graphs.

\begin{thm}\label{25}
Let $G(m_1,\cdots, m_l ;  ~n_1,\cdots,n_l)$ be a chain graph. Then the following hold, $$\sum\limits_{i=1}^{l}m_i+n_1\leq \lambda(G)\leq \sum\limits_{i=1}^{l}m_i+n_1+\sum\limits_{i=1}^{l-1}(n_{i+1}-m_{l+1-i}).$$ If $n_{i+1}>m_{l+1-i}$, then the maximum is attained and if $n_{i+1}\leq m_{l+1-i}$, then the minimum is attained, where $1\leq i\leq {l-1}$. 
\end{thm}
\begin{prf}
Since $G$ is a bipartite graph with partite sets $U$ and $V$, where $U$ is partitioned into disjoint sets $U_1, U_2, \cdots, U_l$ and $V$ is partitioned into $V_1, V_2, \cdots.,V_l$ with cardinalities $m_i$ and $n_i$, respectively for $1\leq i\leq l$. We specify the vertices of $U_1, \cdots, U_l$ and $V_1, \cdots, V_l$ as follows,
\begin{eqnarray*}
U_1&=&\lbrace u_{1_{1}}, u_{2_{1}}, \cdots, u_{m_{1}}\rbrace,\\
U_2&=&\lbrace u_{1_{2}}, u_{2_{2}}, \cdots, u_{m_{2}}\rbrace,\\
\vdots\\
U_l&=&\lbrace u_{1_{l}}, u_{2_{l}}, \cdots, u_{m_{l}}\rbrace,\\
V_1&=&\lbrace v_{1_{1}}, v_{2_{1}}, \cdots, v_{n_{1}}\rbrace,\\
V_2&=&\lbrace v_{1_{2}}, v_{2_{2}}, \cdots, v_{n_{2}}\rbrace,\\
\vdots\\
V_l&=&\lbrace v_{1_{l}}, v_{2_{l}}, \cdots, v_{n_{l}}\rbrace.
\end{eqnarray*}
We define $L(2, 1)$-coloring $f : V(G)\rightarrow [0, \infty)$ of $G$ as,
\begin{eqnarray*}
f(u_{i_{1}}) &=& i_1-1,~~~~~ 1\leq i_1\leq m_1,\\ f(u_{i_{2}}) &=& m_1+i_2, ~~~ 1\leq i_2\leq m_2,\\\vdots\\
f(u_{i_{l}}) &=& \sum\limits_{i=1}^{l-1}(m_i)-1+i_l,~~ 1\leq i_l\leq m_l,\\
f(v_{i_{1}}) &=&\sum\limits_{i=1}^{l}(m_i)+i_1,~~~ 1\leq i_1\leq n_1,\\
f(v_{1_{2}}) &=& \sum\limits_{i=1}^{l}m_i.
\end{eqnarray*}
Note that the vertex $v_{1_{2}}$ is colored by the hole $\sum\limits_{i=1}^{l} m_i$.
The coloring for rest of the vertices of $V_2$ depends on the cardinality of $U_l$, since $d(u_{i_{l}}, v_{i_{2}})=3$ for each $u_{i_{l}}$ and  $v_{i_{2}}$. Similarly, coloring of vertices of $V_3$ depends on the cardinalities of $U_l$ and $U_{l-1}$ and in general, the coloring of vertices of $V_l$ depends on the cardinalities of sets $U_2, \cdots, U_l$.
Suppose, $|V_2|\leq |U_l|$. Then color all but one ($v_{1_{2}}$) vertices of $V_2$ by the existing colors of vertices of $U_l$. In fact, if for each $i$, $3\leq i\leq l$,  $|V_i|\leq |U_{l+2-i}|$, then color the vertices of $V_i$ by the existing colors of $U_{l+2-i}$. Therefore by the definition, $f$ is well-defined minimal coloring. Thus, $span(f)=\sum\limits_{i=1}^{l}(m_i)+n_1$ and consequently, $\lambda(G)=\sum\limits_{i=1}^{l}(m_i)+n_1$.\\
$~~~~$For each $2\leq i\leq l$, suppose $|V_i|> |U_{l+2-i}|$. Then, $|V_2|- |U_{l}|= r_{1}, |V_3|- |U_{l-1}|=r_{2}, \cdots, |V_l|- |U_{2}|=r_{k}$. By Algorithm-(b), we define $L(2, 1)$-colors of $r_1, \cdots, r_k$ vertices as,
   \begin{eqnarray*}
f({v_{{i_{2}}^{'}}})&=&\sum\limits_{i=1}^{l}(m_i)+n_1+{i_{2}^{'}},~~~~1\leq {i_{2}^{'}}\leq r_1,\\
f({v_{{i_{3}}^{'}}})&=&\sum\limits_{i=1}^{l}(m_i)+n_1+r_1+{i_{3}^{'}},~~~~1\leq {i_{3}^{'}}\leq r_2,\\
 \vdots\\
f({v_{{i_{l}}^{'}}})&=&\sum\limits_{i=1}^{l}(m_i)+n_1+r_1+\cdots+r_{k-1}+{i_{l}^{'}},~~~~1\leq {i_{l}^{'}}\leq r_k.
\end{eqnarray*} 
Again, by the definition, $f$ is the well-defined minimal coloring of $G$ and therefore, $\lambda(G)=span(f)=\sum\limits_{i=1}^{l}(m_i)+n_1+r_1+\cdots+r_k=\sum\limits_{i=1}^{l}(m_i)+n_1+(n_2-m_l)+\cdots+(n_l-m_2)=\sum\limits_{i=1}^{l}m_i+n_1+\sum\limits_{i=1}^{l-1}(n_{i+1}-m_{l+1-i})$.
\end{prf}

We conclude this paper with the following example which illustrates Algorithm-(b) and Theorem \eqref{25}.
\begin{exm}
Let $G(5, 2, 3, 2, 4, 3, 3; 2, 2, 5, 3, 6, 2, 5)$ be a chain graph. Define the $L(2, 1)$-coloring $f : V(G)\rightarrow [0, \infty)$ of $G$ as,
\begin{eqnarray*}
f(u_{i_{1}}) &=& i_1-1,~~~~ 1\leq i_1\leq 5,\\  f(u_{i_{2}}) &=& 4+i_2,~~~~ 1\leq i_2\leq 2,\\
f(u_{i_{j}}) &=& 4+\sum\limits_{j=3}^{7}(m_{j-1})+i_j,~~~~ 1\leq i_j\leq m_j~~and~~ 3\leq j\leq 7,\\
f(v_{i_{1}}) &=& \sum\limits_{i=1}^{7}(m_{i})+i_1=22+i_1, ~~~~ 1\leq i_1\leq 2.
\end{eqnarray*}
We compute $|V_2|-|U_7|=-1 = s_1$ (say), it follows that the vertices of $V_2$ can be colored by already existing colors of vertices of $U_7$. Therefore by Algorithm-(b), we define the colors of vertices of $V_2$ as,
\begin{eqnarray*}
f(v_{1_{2}}) &=& \sum\limits_{i=1}^{7}(m_{i})=22,\\
f(v_{2_{2}}) &=& \sum\limits_{i=1}^{7}(m_{i})-1=21.
\end{eqnarray*}
Next, we compute $|V_3|-|U_6|=2=r_1$ and $r_1+s_1=1$. Then by Algorithim-(b), we define the coloring of vertices of $V_3$ as, 
\begin{eqnarray*}
f(v_{i_{3}}) &=& \sum\limits_{i=1}^{7}(m_{i})-(i_3+1)=22,~~~ 1\leq i_3\leq 4\\\\
~~~~~f(v_{5_{3}}) &=& \sum\limits_{i=1}^{7}(m_{i})+n_1+1=25.
\end{eqnarray*}
The coloring for rest of the vertices of $V$ can be defined as,
\begin{eqnarray*}
f(v_{i_{4}}) &=& \sum\limits_{i=1}^{7}(m_{i})-|V_3|-i_4,~~~~ 1\leq i_4\leq 3,\\
f(v_{i_{5}}) &=& \sum\limits_{i=1}^{7}(m_{i})-|V_3|-|V_4|-i_5,~~~~ 1\leq i_5\leq 3,\\
f(v_{i_{5}^{'}}) &=& \sum\limits_{i=1}^{7}(m_{i})+n_1+1+{i_{5}^{'}},~~~~ 1\leq {i_{5}^{'}}\leq 3,\\
\end{eqnarray*}
\begin{eqnarray*}
f(v_{i_{6}}) &=& \sum\limits_{i=1}^{7}(m_{i})-|V_3|-|V_4|- |V_5-3|-i_6,~~~~ 1\leq i_6\leq 2,\\
f(v_{i_{7}}) &=& \sum\limits_{i=1}^{7}(m_{i})-|V_3|-|V_4|- |V_5-3|-|V_6|-i_7,~~~~ 1\leq i_7\leq 3,\\
f(v_{i_{7}^{'}}) &=& \sum\limits_{i=1}^{7}(m_{i})+n_1+1+3+{i_{7}^{'}},~~~~ 1\leq {i_{7}^{'}}\leq 2.
\end{eqnarray*}
Therefore, $\lambda(G)=span(f)=\sum\limits_{i=1}^{7}(m_{i})+n_1+1+3+2=30$. 
\end{exm}

\end{document}